\documentstyle[12pt,leqno,amsmath,amsfonts,amssymb]{article}
\newtheorem{theorem}{Theorem}[section]
\newtheorem{definition}[theorem]{Definition}

\newtheorem{proposition}[theorem]{Proposition}
\newtheorem{remark}[theorem]{Remark}

\newtheorem{lemma}[theorem]{Lemma}
\newcommand {\Mc}      {{\mathcal M}}
\newcommand {\Uc}      {{\mathcal U}}
\newcommand {\Vc}      {{\mathcal V}}
\newcommand {\Wc}      {{\mathcal W}}
\newcommand {\Hc}      {{\mathcal H}}
\newcommand {\Ac}      {{\mathcal A}}
\newcommand {\Tc}      {{\mathcal T}}
\newcommand {\Bc}      {{\mathcal B}}
\newcommand {\Cc}      {{\mathcal C}}

\newcommand {\tK}      {\widetilde{K}}
\newcommand {\bx}      {\hspace{10mm}$\Box$}
\newcommand {\BX}      {\hspace{10mm}\Box}
\newcommand {\tu}      {\tilde{u}}
\newcommand {\tg}      {\tilde{g}}
\newcommand {\tr}      {\tilde{r}}

\newcommand {\ty}      {\tilde{y}}
\newcommand {\ve}      {\varepsilon}
\newcommand {\R}       {{\bf R}}
\newcommand {\RN}      {\R^n}
\newcommand {\MPS}     {M^{1,p}(S,d,\mu)}
\newcommand {\BB}      {\widetilde{B}}
\newcommand {\BT}      {\widetilde{\Bc}_K}
\newcommand {\cw}      {\curlywedge}
\newcommand {\CR}      {C_{rd}}
\newcommand {\br}[1]   {#1^\curlywedge}
\newcommand {\BK}      {\Bc_K}
\newcommand {\MPX}     {M^{1,p}(X,d,\mu)}
\newcommand {\nn}      {\nonumber}
\newcommand {\rf}[1]   {(\ref{#1})}
\newcommand {\card}    {\operatorname{card}}
\newcommand {\supp}    {\operatorname{supp}}
\newcommand {\Ext}     {\operatorname{Ext_S}}
\newcommand {\diam}    {\operatorname{diam}}
\newcommand {\dist}    {\operatorname{dist}}
\newcommand {\cl}      {\operatorname{cl}}

\newcommand {\SECT}[2] {\section*{\centerline{\normalsize
{\bf #1}}} \setcounter{section}{#2}
\setcounter{theorem}{0}\setcounter{equation}{0}}
\newcommand{\be}            {\begin{eqnarray}}
\newcommand{\bel}[1]        {\begin{eqnarray}\label{#1}}
\newcommand{\ee}            {\end{eqnarray}}
\begin{document}
\medskip
\centerline{\large{\bf On extensions of Sobolev functions}}
\vspace*{4mm}
\centerline{\large{\bf defined on regular subsets of metric
measure spaces}}
\vspace*{10mm}
\centerline{By  {\it Pavel Shvartsman}} \vspace* {3 mm}
\centerline {\it Department of Mathematics, Technion -
Israel Institute of Technology,}
\centerline{\it 32000 Haifa, Israel}\vspace*{3 mm}
\centerline{\it e-mail: pshv@tx.technion.ac.il}
\vspace*{12 mm}
\renewcommand{\thefootnote}{ }
\footnotetext[1]{{\it\hspace{-6mm}Math Subject
Classification} 46E35}
\begin{abstract}
{\small
\par We characterize the restrictions of first order
Sobolev functions to regular subsets of a homogeneous
metric space and prove the existence of the corresponding
linear extension operator.}
\end{abstract}
\SECT{1. Main definitions and results.}{1}
\indent
\par Let $(X,d,\mu)$ be a metric space $(X,d)$ equipped
with a Borel measure $\mu$, which is non-negative and outer
regular, and is finite on every bounded subset. In this
paper we describe the restrictions of first order Sobolev
functions to measurable subsets of $X$ which have a certain
regularity property.
\par There are several known ways of defining Sobolev
spaces on abstract metric spaces,
where of course we cannot use
the notion of derivatives. Of particular interest to us,
among these definitions, is the one introduced by
Haj{\l}asz \cite{H1}. But let us first consider a classical
characterization of classical
Sobolev spaces due to Calder\'on. Since it does not use
derivatives, it can lead to yet another way of defining
Sobolev spaces on metric spaces. In \cite{C}
(see also \cite{CS})
Calder\'on characterizes the Sobolev spaces $W^{k,p}(\RN)$
in terms of $L^p$-properties of sharp maximal functions. To
generalize this characterization to the setting of a metric
measure space $(X,d,\mu)$, let $f$ be a locally integrable
real valued function on $X$ and let $\alpha$ be a positive
number. Then the {\it fractional sharp maximal function}
of $f$, is defined by
$$
f^{\sharp}_\alpha(x):=\sup_{r>0}
\frac{r^{-\alpha}}{\mu(B(x,r))}
\int_{B(x,r)}|f-f_{B(x,r)}|\,d\mu.
$$
Here $B(x,r):=\{y\in X:d(y,x)<r\}$ denotes the open ball
centered at $x$ with radius $r$, and, for every Borel set
$A\subset X$ with $\mu(A)<\infty$, $f_A$ denotes
the average value of $f$ over $A$
$$
f_A:=\frac{1}{\mu(A)}\int_A fd\mu.
$$
If $A=\emptyset$, we put $f_A:=0$.
\par In \cite{C} Calder\'{o}n proved
that, for $1<p\le \infty$, the function $u$ is in
$W^{1,p}(\RN)$, if and only if $u$ and $u^{\sharp}_1$
are both in $L^p(\RN)$. This result motivates us to
introduce the space $CW^{1,p}(X,d,\mu)$, which
we will call the {\it Calder\'{o}n-Sobolev space}.
We define it to consist of all
functions $u$ defined on $X$ such that
$u,u^{\sharp}_1\in L^p(X)$.
We equip this space with the Banach norm
$$
\|u\|_{CW^{1,p}(X,d,\mu)}:=
\|u\|_{L^p(X)}+\|u^{\sharp}_1\|_{L^p(X)}.
$$
\par Let us now recall the details of the definition of
Haj\l asz mentioned above.
Haj\l asz \cite{H1} introduced the Sobolev-type space
on a metric space, $M^{1,p}(X,d,\mu)$ for $1<p\le\infty$.
It consists of all functions $u\in L^p(X)$ for which there
exists a function $g\in L^p(X)$ (depending on $u$)
such that the inequality
\bel{Pin}
|u(x) - u(y)|\le d(x, y)(g(x) + g(y))
\ee
holds $\mu$-a.e. (This means that  there is a set
$E\subset X$ with $\mu(E)=0$ such that \rf{Pin}
holds for every $x,y\in X\setminus E$).
As in \cite{HKin} we will refer to all functions $g$
which satisfy the inequality \rf{Pin} as
{\it generalized gradients} of $u$.
\par $M^{1,p}(X,d,\mu)$ is normed by
$$
\|u\|_{M^{1,p}(X,d,\mu)}:=\|u\|_{L^p(X)}
+\inf_{g}\|g\|_{L^p(X)}
$$
where the infimum is taken over
the family of all generalized gradients of $u$.
\par In the case where $X=\Omega\subset\RN$ is an
open bounded domain with a Lipschitz boundary, $d$ is the
Euclidean distance and $\mu$ is the $n$-dimensional
Lebesgue measure on $\Omega$, Haj\l asz \cite{H1} showed
that the space $M^{1,p}(\Omega,d,\mu)$ coincides with the
Sobolev space $W^{1,p}(\Omega)$ and, moreover, that every
function $u\in W^{1,p}(\RN)$ satisfies \rf{Pin} with
$g=c\Mc\|\nabla u\|$. Here $\Mc$ is the Hardy-Littlewood
maximal operator and $c=c(n)$. (For further development and
application of this approach to Sobolev spaces on metric
space see, e.g. \cite{FHK,FLW,GT,H1,HKin,HKos,HM,He,KL} and
references therein.)
\par It turns out that for a {\it doubling} measure $\mu$,
the Haj{\l}asz-Sobolev space coincides with the
Calder\'on-Sobolev space, i.e.,
$$
CW^{1,p}(X,d,\mu)=M^{1,p}(X,d,\mu), ~~~1<p\le\infty,
$$
and, moreover, for every $u\in M^{1,p}(X,d,\mu)$,
the function $g=cu^{\sharp}_1$
(with some constant $c=c(X)$) is a generalized
gradient of $u$. This is an immediate consequence of
a result of Haj\l asz and Kinnunen.
(See \cite{HKin}, Theorem 3.4.)
\par We recall that a measure $\mu$ satisfies the
{\it doubling condition} if there exists
a constant $C_d\ge 1$ such that, for every
$x\in X$ and $r>0$,
\bel{double}
\mu(B(x,2r))\le C_d\mu(B(x,r)).
\ee
\par As usual, see \cite{CW},
we will call a metric measure space $(X,d,\mu)$
with a doubling measure
$\mu$ a {\it metric space of homogeneous type}
and $C_d$ a {\it doubling constant}.
\par In this paper we will only consider such metric
measure spaces, which also satisfy
an additional condition, namely that there exists
a constant $C_{rd}>1$
such that, for every $x\in X$ and $r>0$,
\bel{mreg}
C_{rd}\mu(B(x,r))\le \mu(B(x,2r)).
\ee
We call this condition the reverse doubling condition and
$C_{rd}$ a reverse doubling constant.
\par We will characterize the restrictions
of Calder\'{o}n-Sobolev and Haj{\l}asz-Sobo\-lev functions
to {\it regular} subsets of a homogeneous metric space
$(X,d,\mu)$.
\begin{definition}\label{REG}
A measurable set $S\subset X$ is said to be regular if
there are constants $\theta_S\ge 1$ and $\delta_S>0$ such
that for every $x\in S$ and $0<r\le\delta_S$
$$
\mu(B(x,r))\le \theta_S \mu(B(x,r)\cap S).
$$
\end{definition}
\par  A Cantor-like set or a Sierpi\'{n}ski's type gasket
(or carpet) of positive Lebes\-gue measure provide examples
of non-trivial regular subsets of $\RN$. (Regular subsets
of $\RN$ are often called Ahlfors $n$-regular or $n$-sets
\cite{JW}). For properties of metric spaces supporting
doubling measures and sets satisfying regularity conditions
we refer to \cite{BG,J,JW,VK} and references therein.
\par Let us formulate the main result of the paper.
Given a Borel set $A\subset X$, a function $f\in
L^{1,loc}(A)$ and $\alpha>0$ we let $f^{\sharp}_{\alpha,A}$
denote the fractional sharp maximal function of $f$ on $A$,
\bel{smfS} f^{\sharp}_{\alpha,A}(x):=\sup_{r>0}
\frac{r^{-\alpha}}{\mu(B(x,r))} \int_{B(x,r)\cap A}
|f-f_{B(x,r)\cap A}|\,d\mu,~~~~x\in A.
\ee
Thus, $f^{\sharp}_{\alpha}=f^{\sharp}_{\alpha,X}$~.
\par As usual for a Banach space $(\Ac, \|\cdot\|_{\Ac})$
of measurable functions defined on $X$ and a Borel
set $S\subset X$ we let $\Ac|_S$ denote the
restriction of $\Ac$ to $S$, i.e., a Banach space
$$
\Ac|_S:=\{f:S\to\R~: {\rm there~is}~ F\in \Ac
~{\rm such~that}~ F|_S=f\}
$$
equipped with the standard quotient space norm
$$
\|f\|_{\Ac|_S}:=\inf\{\|F\|_{\Ac}: F\in \Ac,F|_S=f\}.
$$
\begin{theorem}\label{EXT1}
Let $(X,d,\mu)$ be a metric space of homogeneous type
satisfying the reverse doubling condition \rf{mreg}
and let $S$ be a regular subset of $X$.
Then a function $u\in L^p(S),~1<p\le \infty,$
can be extended to a function
$\tilde{u}\in CW^{1,p}(X,d,\mu)$ if and only if
$u_{1,S}^{\sharp}\in L^p(S)$.
In addition,
$$
\|u\|_{CW^{1,p}(X,d,\mu)|_S}\approx
\|u\|_{L^p(S)}+\|u_{1,S}^{\sharp}\|_{L^p(S)}
$$
with constants of equivalence depending
only on $C_d,\CR,\theta_S,\delta_S$ and $p$.
Moreover, there exist a linear continuous
extension operator
$$
\Ext:CW^{1,p}(X,d,\mu)|_S\to CW^{1,p}(X,d,\mu).
$$
Its operator norm is bounded by a constant depending
only on $C_d,\CR,\theta_S,\delta_S$ and $p$.
\end{theorem}
\par Let us apply this result to $X=\RN$ with the Lebesgue
measure (clearly, in this case \rf{mreg} is satisfied with
$\CR=2^n$). Then for every regular subset $S\subset\RN$ we
have:
\par (i).~~
$W^{1,p}(\RN)|_{S}= \{u:S\to\R~:~u,u_{1,S}^{\sharp}\in
L^p(S)\}, ~~~~1<p\le\infty$;\vspace{2mm}
\par (ii). There is a linear continuous
extension operator from
$W^{1,p}(\RN)|_{S}$ into $W^{1,p}(\RN)$.\vspace{2mm}
\par Observe that (ii) follows from a general
result of Rychkov \cite{R}.
\par There is an extensive literature devoted to
description of the restrictions of Sobolev functions to
different classes of subsets of $\RN$. We refer the reader
to the books of Maz'ya \cite{M}, Maz'ya and Poborchi
\cite{MP}, the article of Farkas and Jakob \cite{FJ} and
references therein for numerous results and technique in
this direction. We also observe that the criterion (i) can
be useful for description of Sobolev extension domains,
i.e., domains $\Omega\subset \RN$ such that
$W^{1,p}(\RN)|_\Omega=W^{1,p}(\Omega)$. For instance, due
to a result of Koskela \cite{Kos}, every Sobolev extension
domain is a regular subset of $\RN$ whenever
$n-1<p<\infty$.
\par The second main result of the paper is the following
\begin{theorem}\label{EXT2}
Let $(X,d,\mu)$ be a homogeneous
metric space satisfying condition \rf{mreg}.
Then for every regular subset $S$ of $X$
$$
M^{1,p}(X,d,\mu)|_S=M^{1,p}(S,d,\mu).
$$
Moreover, there exist a linear continuous extension
operator
$$
\Ext:~M^{1,p}(S,d,\mu)\to M^{1,p}(X,d,\mu)
$$
such that
$\|\Ext\|\le C(C_d,\CR,\theta_S,\delta_S,p)$.
\end{theorem}
\par For a family of bounded domains in $\RN$ satisfying
a certain plumpness condition (so-called $A(c)$-condition)
Theorem \ref{EXT2} was proved by Haj\l asz and Martio
\cite{HM}. Harjulehto \cite{Har} has generalized this
result for the case of homogeneous metric spaces
$(X,d,\mu)$ and domains $\Omega$ satisfying so-called
$A^*(\varepsilon,\delta)$-condition. Observe that both
$A(c)$- and $A^*(\varepsilon,\delta)$-sets are regular but
a Cantor-type set of positive Lebesgue measure in $\RN$
provides an example of a regular subset which satisfies
neither $A(c)$- nor $A^*(\varepsilon,\delta)$-condition.
\par Proofs of Theorems \ref{EXT1} and \ref{EXT2}
are based on a modification of the Whitney extension method
suggested in author's work \cite{S}, see also \cite{S2},
for the case of regular subsets of $\RN$.
\par A crucial step of this approach is presented in
Section 2. Without loss of generality we may assume that
$S$ is closed (see Lemma \ref{DENS}) so that $X\setminus S$
is open. Since $\mu$ is doubling, $X\setminus S$ admits a
Whitney covering which we denote by $W_S$.
\par We assign every ball $B=B(x_B,r_B)\in W_S$
a measurable subset $H_B\subset S$ such that $H_B\subset
B(x_B,\gamma_1r_B)\cap S$, $\mu(B)\le \gamma_2\mu(H_B)$
whenever $r_B\le \delta_S$, and the family
$$
\Hc_S:=\{H_B:~B\in W_S\}
$$
has {\it a finite covering multiplicity }, i.e., every
point $x\in S$ belongs at most $\gamma_3$ sets of the
family $\Hc_S$. Here $\gamma_1,\gamma_2,\gamma_3$ are
positive constants depending only on $C_d,\CR$ and
$\theta_S$. We call every set $H_B\in\Hc_S$ a {\it
``reflected quasi-ball"} associated to the Whitney ball
$B$. The existence of the family $\Hc_S$ of reflected
quasi-balls is proved in Theorem \ref{HB}.
\par The second step of the extension method and the proof
of Theorem \ref{EXT2} are presented in Section 3. We fix
functions $u\in M^{1,p}(S,d,\mu)$ and $g\in L^p(S)$
satisfying on $S$ inequality \rf{Pin}. Then we define an
extension $\tu$ of $u$ by the formula
\bel{OE}
\tilde{u}(x)=(\Ext u)(x):=\sum_{B\in W_S}
u_{H_B}\varphi_B(x), ~~~x\in X\setminus S.
\ee
Here $\{\varphi_B:B\in W_S\}$ is a partition of unity
associated to the Whitney covering.
\par Finally, we define an extension $\tg$ of $g$
by letting
$$
\tg(x):=\sum_{B\in W_S}\limits (g_{H_B}+|u_{H_B}|)
\chi_{B^*}(x),
~~~x\in X\setminus S,
$$
where $B^*:=B(x_B,\frac{9}{8}r_B)$. We show that
$\tilde{u}\in L^p(X)$, $\tilde{g}\in L^p(X)$
and $\tilde{g}$ is a generalized gradient of
$\tilde{u}$, i.e., the pair $(\tilde{u},\tilde{g})$
satisfies on $X$ inequality \rf{Pin}. Since
$\tilde{u}|_S=u$, this proves that
$u\in M^{1,p}(X,d,\mu)|_S$ so that $\Ext$ provides a
{\it linear extension operator} from $\MPS$ into $\MPX$.
\par Section 4 is devoted to estimates of the sharp maximal
function of the extension $\tilde{u}:=\Ext u$. Given a
function $f$ defined on $S$ we let $f^\curlywedge$ denote
its extension on all of $X$ by zero. We show that for every
$\alpha>0$ and $x\in X$
$$
(\tu)^\sharp_\alpha(x)\le
C(\Mc(u^\sharp_{\alpha,S})^\cw(x)+
\Mc u^\cw(x)),
$$
see Theorem \ref{pr}. Using this estimate and the
Hardy-Littlewood maximal theorem we prove a slightly more
general version of Theorem \ref{EXT1} related to the space
$\Cc^\alpha_p(X,d,\mu)$. This space consists of all
functions $u$ defined on $X$ such that $u,
u^\sharp_\alpha\in L^p(X)$. $\Cc^\alpha_p(X,d,\mu)$ is
normed by
$$
\|u\|_{\Cc^\alpha_p(X,d,\mu)}:=
\|u\|_{L^p(X)}+\|u^{\sharp}_\alpha\|_{L^p(X)}.
$$
For the case $X=\RN$ with the Lebesgue measure this
space was introduced and investigated by
Devore and Sharpley \cite{DS} and Christ \cite{Ch}.
Clearly, $CW^{1,p}(X,d,\mu)=\Cc^1_p(X,d,\mu)$.
\begin{theorem}\label{EXT3}
Let $(X,d,\mu)$ be a metric space of homogeneous type
satisfying condition \rf{mreg} and let $S$ be a regular
subset of $X$. A function $u\in L^p(S),~1<p\le \infty,$
belongs to the trace space $\Cc^\alpha_p(X,d,\mu)|_S$ if
and only if $u_{\alpha,S}^{\sharp}\in L^p(S)$. In addition,
\bel{Eq}
\|u\|_{\Cc^\alpha_p(X,d,\mu)|_S}\approx
\|u\|_{L^p(S)}+\|u_{\alpha,S}^{\sharp}\|_{L^p(S)}
\ee
and there exists a linear continuous extension operator
$$
\Ext:~\Cc^\alpha_p(X,d,\mu)|_S
\to \Cc^\alpha_p(X,d,\mu)
$$
whose operator norm is bounded by a constant depending only
on $C_d,\CR,\theta_S,\delta_S$ and $p$.
\end{theorem}
\par Observe that for $X=\RN$ and $S$ to be a Lipschitz or
an $(\varepsilon,\delta)$-domain this result follows from
extension theorems proved by
Devore and Sharpley \cite{DS}, pp. 99--101,
(Lipschitz domains), and Christ \cite{Ch}
($(\varepsilon,\delta)$-domains).
\par {\bf Acknowledgement.} I am very grateful to
M. Cwikel for helpful discussions and valuable advice.
\SECT{2. The Whitney covering and a family of associated
``quasi-balls".}{2}
\indent
\par We will use the following notation.
Throughout the paper $C,C_1,C_2,...$ will be generic
positive constants which depend only on
$C_d,\CR,\theta_S,\delta_S$ and $p$. These constants can
change even in a single string of estimates. We write
$A\approx B$ if there is a constant $C$ such that $A/C\le
B\le CA$. For a ball $B=B(x,r)$ we let $x_B$ and $r_B$
denote center and radius of $B$. Given a constant $\lambda>
0$ we let $\lambda B$ denote the ball $B(x,\lambda r)$. For
$A,B\subset X$ and $x\in X$ we put
$$
\dist(A,B):=\inf\{d(a,b):~a\in A, b\in B\}
$$
and $d(x,A):=\dist(\{x\},A)$. Finally, by $\cl(A)$ we
denote the closure of $A$ in $X$.
\begin{lemma}\label{DENS}
For every regular subset $S\subset X$
$$\mu(\cl(S)\setminus S)=0.$$
\end{lemma}
{\bf Proof.} Denote $Y:=\cl(S)\setminus S$ and fix $y\in
Y$. Then for every $r,~0<r\le\delta,$ there is a point
$\ty\in S$ such that $\dist(y,\ty)\le r/4$. Clearly,
$B(\ty,r/4)\subset B(y,r)$. Since $S$ is regular and
$\ty\in S$, we obtain
$$
\mu(B(y,r)\cap S)\ge\mu(B(\ty,r/4)\cap S)
\ge\theta\mu(B(\ty,r/4)).
$$
On the other hand,
$B(y,r)\subset B(\ty,5r/4)$ so that
by the doubling condition
$$
\mu(B(y,r))\le\mu(B(\ty,5r/4))
\le C_d^3\mu(B(\ty,r/4)).
$$
Hence $\mu(B(y,r)\cap S)\ge\theta C_d^{-3}\mu(B(y,r))$.
We let $D_A$ denote the family of density
points of the set $A:=X\setminus S$. Then
$$
\frac{\mu(B(y,r)\cap A)}{\mu(B(y,r))}
<1-\theta C_d^{-3},~~~~~ y\in Y,
$$
which implies $Y\cap D_A=\emptyset$. Thus
$Y\subset A\setminus D_A$ so that by Lebesgue's theorem,
see, e.g. \cite{F}, \S 2.9,
$\mu(Y)\le\mu(A\setminus D_A)=0$.\bx
\par In the remaining part of the paper we will assume
that $S$ is a {\it closed} regular subset of $X$.
\par Since $\mu$ is a doubling measure,
there exists a constant $M=M(C_d)$ such that in every ball
$B(x,r)$ there are at most $M$ points ${x_j}$ satisfying
the inequality $d(x_i,x_j)\ge \frac{r}{2}$ (one can put
$M:=C^4_d$). For every metric space with this property the
following is true (see, e.g. \cite{G}, Theorem 2.3): For
every open subset $G\subset X$ with a non-empty boundary
there is a countable family of balls $\Wc_G$ such that
$G=\cup\{B:B\in\Wc_G\}$, every point of $G$ is covered by
at most $9M$ sets from $\Wc_G$ and $r\le
\dist(B(x,r),\partial G)\le 4r$ for every $B=B(x,r)\in
\Wc_G$.
\par Applying this result to the open set
$G=X\setminus S$ we obtain the following
\begin{theorem}\label{Wcov}
There is a countable family of balls $W_S$ such that
\par (i). $X\setminus S=\cup\{B:B\in W_S\}$;
\par (ii). For every ball $B=B(x,r)\in W_S$
$$
r\le \dist(B(x,r),S)\le 4r;
$$
\par (iii). Every point of $X\setminus S$ is covered by
at most $N=N(C_d)$ balls from $W_S$.
\end{theorem}
\par Using standard argument one can readily prove
the following additional properties of Whitney's balls.
\begin{lemma}\label{Wadd}
\par (1). For every $B=B(x_B,r_B)\in W_S$
there is a point $y_B\in S$ such that
\bel{inBY}
B(y_B,r_B)\subset B(x_B,Cr_B) {\rm~~~and~~~}
B=B(x_B,r_B)\subset B(y_B,Cr_B).
\ee
Moreover, $\mu(B(x_B,r_B))\approx\mu(B(y_B,r_B))$.
\par (2). If $B,K\in W_S$ and $B^*\cap K^*\ne\emptyset$,
then
\bel{eqvball}
C^{-1}r_B\le r_K\le Cr_B.
\ee
(Recall that $B^*:=\frac{9}{8}B$.)
\par (3). For every ball $K\in W_S$ there are at most
$N$ balls from the family $W^*_S:=\{B^*:B\in W_S\}$ which
intersect $K^*$.
\par (4). For every $B\in W_S$ and every $x\in B^*$
\bel{drx}
C^{-1}r_B\le d(x,S)\le C r_B.
\ee
\par Here $C$ and $N$ are positive constants
which depend only on $C_d$.
\end{lemma}
\par The next lemma easily
follows from inequalities \rf{double} and \rf{mreg}.
\begin{lemma}\label{CDCR} For every $x\in X, r>0,
1\le t<\infty,$
\bel{muCr}
\mu(B(x,r))\le \CR t^{-\alpha}\mu(B(x,tr))
\ee
and
\bel{muCd}
\mu(B(x,tr))\le C_dt^\beta\mu(B(x,r))
\ee
where $\alpha:=\log_2\CR$ and $\beta:=\log_2C_d.$
\end{lemma}
\par Let us formulate the main result of the section.
\begin{theorem}\label{HB} There is a family of Borel
sets $\Hc_S=\{H_B:~B\in W_S\}$ such that:
\par (i). $H_B\subset (\gamma_1 B)\cap S,~~~~~B\in W_S$;
\par (ii). $\mu(B)\le \gamma_2\mu (H_B)$ whenever
$B\in W_S$ and $r_B\le\delta_S$;
\par (iii). $\sum_{B\in W_S}\limits\chi_{H_B}\le \gamma_3$.
\par Here  $\gamma_1,\gamma_2,\gamma_3$ are positive
constants depending only on $\CR,C_d$ and $\theta_S$.
\end{theorem}
{\bf Proof.} Let $K=B(x_K,r_K)\in W_S$ and let $y_K$ be a
point on $S$ satisfying condition (1) of Lemma \ref{Wadd}.
Thus $B(y_K,r_K)\subset CK$ and $K\subset B(y_K,Cr_K)$.
\par Given $\ve, 0<\ve\le 1,$ we denote
$K_{\ve}:=B(y_K,\ve r_K)$. Let $B=B(x_B,r_B)$ be a ball
from $W_S$ with
$r_B\le \delta_S$. Set
\bel{2.13'}
\Ac_B:=\{K=B(x_K,r_K)\in W_S:~K_\ve\cap
B_\ve\ne\emptyset,~r_K\le \ve r_B\}.
\ee
Recall that $B_{\ve}:=B(y_B,\ve r_B)$. We define a
``quasi-ball" $H_B$ by letting
\bel{dAB}
H_B:=(B_\ve\cap S)\setminus
(\cup\{K_\ve:~K\in \Ac_B\}).
\ee
If $r_B>\delta_S$ we put $H_B:=\emptyset$.
\par Prove that for some $\ve:=\ve(\CR,C_d,\theta)$
small enough the family of subsets $\Hc_S:=\{H_B:~B\in
W_S\}$ satisfies conditions (i)-(iii). By \rf{dAB} and
\rf{inBY}
$$
H_B\subset B_\ve:=B(y_B,\ve r_B)
\subset B(y_B,r_B)\subset B(x_B,Cr_B)=CB.
$$
In addition, by \rf{dAB} $H_B\subset S$ so that $H_B\subset
(CB)\cap S$ proving property (i).
\par Let us prove (ii). Suppose that $B=B(x_B,r_B)\in W_S$
and $r_B\le\delta_S$.
If $K\in \Ac_B$, then by \rf{2.13'}
$K_\ve\cap B_\ve\ne \emptyset$ and $r_K\le\ve r_B$. Hence
$$
r_{K_\ve}(=\ve r_K)\le\ve r_{B_\ve}(:=\ve^2r_B)
\le r_{B_\ve}
$$
so that $y_K\in 2B_\ve$. But
$r_K\le\ve r_B=r_{B_\ve}$ and $K\subset B(y_K,Cr_K)$ which
implies $K\subset(C+2)B_\ve$. Thus
\bel{2.**}
\Uc_B:=\cup\{K:~K\in \Ac_B\}\subset(C+2)B_\ve.
\ee
By property (iii) of Theorem \ref{Wcov}
$$
\sum_{K\in \Ac_B}\limits\chi_K(x)\le
\sum_{K\in W_S}\limits\chi_K(x)\le
N=N(C_d),~~~~~x\in X,
$$
so that by \rf{2.**} and \rf{muCd}
$$
\sum_{K\in \Ac_B}\mu(K)=
\int_{\Uc_B}
\sum_{K\in \Ac_B}\chi_K\, d\mu
\le\int_{(C+2)B_\ve} Nd\mu=N\mu((C+2)B_\ve)
\le C_1\mu(B_\ve).
$$
On the other hand, for every $K\in \Ac_B$ by \rf{muCr}
and by (1), Lemma \ref{Wadd}
$$
\mu(K_\ve)=\mu(B(y_K,\ve r_K))\le
\CR\ve^{\alpha}\mu(B(y_K,r_K))
\le C_2\ve^{\alpha}\mu(K).
$$
Hence
$$
\mu(\cup\{K_\ve:~K\in \Ac_B\})
\le
\sum_{K\in \Ac_B}\mu(K_\ve)\\
\le
C_2\ve^{\alpha}\sum_{K\in \Ac_B}\mu(K)
\le
C_3\ve^{\alpha}\mu(B_\ve).
$$
\par Since $S$ is regular and
$r_{B_\ve}=\ve r_B\le \delta_S$,
$\mu(B_\ve\cap S)\ge \theta^{-1} \mu(B_\ve)$ so that
\be
\mu(H_B)&=&\mu((B_\ve\cap S)
\setminus (\cup\{K_\ve:~K\in \Ac_B\}))\nn\\
&\ge&
\mu(B_\ve\cap S)-
\mu(\cup\{K_\ve:~K\in \Ac_B\})\ge
(\theta-C_3\ve^{\alpha})\mu(B_\ve).\nn
\ee
By \rf{muCd} and by property (1) of Lemma \ref{Wadd}
\be
\mu(B_\ve)=\mu(B(y_B,\ve r_B))
&\ge&
C_d^{-1}\ve^\beta\mu(B(y_B,r_B))\nn\\
&\ge&
C^{-1}C_d^{-1}\ve^\beta\mu(B(x_B,r_B))=
C_4\ve^\beta\mu(B)\nn
\ee
so that
$$
\mu(H_B)\ge
C_4(\theta^{-1}-C_3\ve^{\alpha})\ve^\beta\mu(B).
$$
We define $\ve$ by setting
$\ve:=(2C_3\theta)^{-\frac{1}{\alpha}}$.
Then the inequality $\mu(B)\le\gamma_2\mu(H_B)$ holds with
$\gamma_2:=2C_4^{-1}\theta^{\frac{\beta}{\alpha}-1}
(2C_3)^{\frac{\beta}{\alpha}}$
proving property (ii) of the theorem.
\par Let us prove (iii). Let
$B=B(x_B,r_B),B'=B(x_{B'},r_{B'})\in W_S$
be Whitney's balls such that $r_B,r_{B'}\le\delta_S$ and
$H_B\cap H_{B'}\ne \emptyset$. Since $H_B\subset B_\ve,
H_{B'}\subset B'_\ve$, we have $B_\ve\cap B'_\ve\ne
\emptyset$.
\par On the other hand,
$B\notin \Ac_{B'}$ and $B'\notin \Ac_{B}$, otherwise by
\rf{2.13'} and \rf{dAB} $H_B\cap H_{B'}=\emptyset$. Since
$B_\ve\cap B'_\ve\ne \emptyset$, by definition \rf{2.13'}
$r_B>\ve r_{B'}$ and $r_{B'}>\ve r_B$ so that $r_B\approx
r_{B'}$. By \rf{inBY}
$$
B_\ve=B(y_B,\ve r_B)\subset B(y_B,r_B)\subset CB
$$
and similarly $B'_\ve\subset CB'$. But $B_\ve\cap B'_\ve\ne
\emptyset$ so that $CB\cap CB'\ne\emptyset$ as well.
Moreover, since $r_B\approx r_{B'}$, we have $B'\subset
C_5B$ and $B\subset C_5B'$. This and the doubling condition
imply $\mu(B')\approx \mu(B)$.
\par We denote
$$
\Tc_B:=\{B'\in W_S:~H_B\cap H_{B'}
\ne\emptyset,~r_{B'}\le\delta_S\}
$$
and $\Vc_B:=\cup\{B':~B'\in\Tc_B\}$.
Thus we have proved that
$\Vc_B\subset C_5B$ and
$\mu(B')\approx \mu(B)$ for every $B'\in\Tc_B$.
\par Let $M_B:=\card\Tc_B$ be the cardinality of $\Tc_B$.
Clearly, to prove (iii) it suffices to show
that $M_B\le\gamma_3$. We have
$$
M_B\mu(B)\le
C\sum_{B'\in\Tc_B}\mu(B')=
C\int_{\Vc_B}\limits
\sum_{B'\in\Tc_B}\chi_{B'}\,d\mu
\le C\int_{C_5B}\limits
\sum_{B'\in\Tc_B}\chi_{B'}\,d\mu.
$$
By the property (iii) of Theorem \ref{Wcov}
$$
\sum\{\chi_{B'}:B'\in\Tc_B\}\le \sum\{\chi_{B'}:B'\in
W_S\}\le N=N(C_d)
$$
so that
$$
M_B\mu(B)\le
C\int_{C_5B}\limits
Nd\mu=C N\mu(C_5B)\le C\mu(B)
$$
proving the required inequality
$M_B\le \gamma_3$.\bx\vspace*{2mm}
\SECT{3. The extension operator: proof of Theorem
\ref{EXT2}.}{3}
\indent
\par For every $u\in \MPX$ and every
generalized gradient $g$ of $u$ the restriction $g|_S$ is
a generalized gradient
of $u|_S$ so that $\MPX|_S\subset \MPS$.
\par Let us prove that formula \rf{OE} provides a linear
continuous extension operator from $\MPS$ into $\MPX$.
Obviously, this will imply the converse imbedding as well.
\par Recall that for every $u\in\MPS$ its generalized
gradient $g$ belongs to $L^p(S)$ and satisfies the
inequality
\bel{lipug}
|u(x) - u(y)|\le d(x, y)(g(x) + g(y)),~~~~~
x,y\in S\setminus E,
\ee
where $E$ is a subset of $S$ of measure $0$. We may suppose
that $g$ is almost optimal, i.e., $\|g\|_{L^p(S)}\le
2\|u\|_{\MPS}$.
\par The extension operator $\Ext$, see \rf{OE},
is determined by the family of Borel subsets
$\Hc_S=\{H_B:~B\in W_S\}$  introduced in the previous
section. We recall that $\mu(H_B)>0$ for every ball $B\in
W_S$ with $r_B\le \delta_S$ and $H_B:=\emptyset$ whenever
$r_B>\delta_S$. Therefore according to our notation
$u_{H_B}$ is the average of $u$ over $H_B$ whenever $r_B\le
\delta_S$ and $u_{H_B}:=0$ otherwise.
\par We let $\Phi_S=\{\varphi_B:~B\in W_S\}$
denote a partition of unity associated to the Whitney
covering $W_S$, see, e.g. \cite{MS}. We recall that
$\Phi_S$ is a family of functions defined on $X$ which have
the following properties: For every ball $B\in W_S$
\par (a). $0\le\varphi_B\le 1$;
(b). $\supp \varphi_B\subset B^*(:=\frac{9}{8}B)$;
(c). $\sum\{\varphi_B(x):~B\in W_S\}=1$ on
$X\setminus S$;
(d). for some constant $C=C(C_d)$
$$
|\varphi_B(x)-\varphi_B(y)| \le C\frac{d(x,y)}{r_B}\,,
~~~~~~x,y\in X.
$$
\par Recall that the extension operator $\tu=\Ext u$
is defined by the formula
\bel{Extu}
\tu(x):=\sum_{B\in W_S}\limits u_{H_B}\varphi_B(x),
~~~~~~x\in X\setminus S,
\ee
and $\tu(x):=u(x)$, $x\in S$. We also define  an extension
$\tg$ of $g$ by letting
\bel{defG}
\tg(x):=
\sum_{B\in W_S}\limits (g_{H_B}+|u_{H_B}|)\chi_{B^*}(x),
~~~~~~x\in X\setminus S,
\ee
and $\tg(x):=g(x)$ for $x\in S$.
\par To prove that $\Ext$ satisfies conditions of
Theorem \ref{EXT2} it suffices to show that
\bel{IN1}
\|\tu\|_{L^p(X)}\le C\|u\|_{L^p(S)},~~~~~~
\|\tg\|_{L^p(X)}\le C(\|g\|_{L^p(S)}+\|u\|_{L^p(S)})
\ee
and the inequality
\bel{IN2}
|\tu(x) - \tu(y)|\le Cd(x, y)(\tg(x) + \tg(y))
\ee
holds $\mu$-a.e. on $X$. Then
$$
\|\tu\|_{\MPX}
\le\|\tu\|_{L^p(X)}+\|\tg\|_{L^p(X)}
\le C(\|u\|_{L^p(S)}+\|g\|_{L^p(S)})
$$
proving that $\|\tu\|_{\MPX}\le C\|u\|_{\MPS}$ and
$\|\Ext\|\le C$.
\par Proofs of inequalities \rf{IN1} and \rf{IN2} are based
on a series of auxiliary lemmas.
\begin{lemma}\label{mean}
Let $H,H'\subset S$ and let $0<\mu(H),\mu(H')<\infty$. Then
\bel{forall}
|u_H - u_{H'}|\le
\diam(H\cup H')(g_{H}+g_{H'})
\ee
and for every $y\in S$
\bel{fory}
|u_{H} - u(y)|\le
\diam(H\cup \{y\})({g}_{H}+g(y)).
\ee
\end{lemma}
{\bf Proof.} We have
$$
I:=|u_H - u_{H'}| \le\frac{1}{\mu(H)}\frac{1}{\mu(H')}
\int_{H}\int_{H'} |u(x)-u(y)|\,d\mu(x)d\mu(y)
$$
so that by \rf{lipug}
$$
I\le \frac{1}{\mu(H)}\frac{1}{\mu(H')}
\int_{H}\int_{H'}
d(x,y)(g(x)+g(y))d\mu(x)d\mu(y).
$$
Since $d(x,y)\le \diam(H\cup H')$ for every
$x\in H,y\in H'$, we have
$$
I\le
\frac{\diam(H\cup H')}{\mu(H)\mu(H')}
\int_{H}\int_{H'}
(g(x)+g(y))d\mu(x)d\mu(y)\nn\\
=
\diam(H\cup H')(g_{H}+g_{H'})
$$
proving \rf{forall}. In a similar way we prove inequality
\rf{fory}.\bx
\begin{lemma}\label{uBB}
Let $\BB\in W_S$ and let $x\in\BB$. Then for every
$y\in X\setminus S$ and every ball $B\in W_S$ such that
$B^*\cap \{x,y\}\ne \emptyset$ we have
\bel{difub}
|u_{H_B}-u_{H_{\BB}}|\le C(d(x,S)+d(x,y)+d(y,S))
(\tg(x)+\tg(y)).
\ee
If $y\in S$, then for every
$B\in W_S$ such that $B^*\ni x$
\bel{difuy}
|u_{H_B}-u(y)|\le Cd(x,y)(\tg(x)+\tg(y)).
\ee
\end{lemma}
{\bf  Proof.} First we prove \rf{difub}.
Suppose that $y\in X\setminus S$ and consider the case
$r_B\le \delta_S, r_{\BB}\le \delta_S$.
Since $\mu(H_B),\mu(H_{\BB})>0$, by \rf{forall}
\be\label{uB1}
|u_{H_B}-u_{H_{\BB}}|\le \diam(H_B\cup H_{\BB})
(g_{H_B}+g_{H_{\BB}}).
\ee
\par By \rf{drx} $r_B\approx d(x,S)$
whenever $x\in B^*$ and
by property (i) of Theorem \ref{HB},
$H_B\subset \gamma_1B$ so that
$H_B\subset B(x,Cd(x,S))$.
Since $x\in\BB\subset (\BB)^*$, we also have
$H_{\BB}\subset B(x,C_2d(x,S))$. In a similar way we prove
that $H_B\subset B(y,C_2d(y,S))$ whenever $y\in B^*$. Hence
$$
\diam(H_B\cup H_{\BB})\le
\diam(B(x,C_2d(x,S))\cup B(y,C_2d(x,S)))
$$
so that
\bel{Dxy}
\diam(H_B\cup H_{\BB})\le
C_2(d(x,S)+d(x,y)+d(y,S)).
\ee
Since $r_B, r_{\BB}\le \delta_S$ and $x\in B^*$ or
$y\in B^*$, by definition of $\tg$, see \rf{defG},
we have $g_{H_{\BB}}\le \tg(x)$,
$g_{H_B}\le \tg(x)$ (whenever $x\in B^*$) or
$g_{H_B}\le \tg(y)$ (if  $y\in B^*$).
Hence
$$
g_{H_B}+g_{H_{\BB}}\le 2(\tg(x)+\tg(y)).
$$
Combining this inequality with \rf{uB1} and \rf{Dxy}
we obtain \rf{difub} for
the case $r_B, r_{\BB}\le \delta_S$.
\par Let us prove \rf{difub} for the case
$r_B>\delta_S,r_{\BB}\le \delta_S$.
By \rf{drx} $\delta_S\le r_B\le Cd(y,S)$
whenever $y\in B^*$ or $\delta_S\le r_B\le Cd(x,S)$,
if $x\in B^*$. Hence
$$
\frac{C}{\delta_S}(d(x,S)+d(x,y)+d(y,S))\ge 1.
$$
Since $x\in \BB$, by \rf{defG}
$|u_{H_{\BB}}|\le \tg(x)$, and since
$r_B>\delta_S$,  $u_{H_B}:=0$. Hence
$$
|u_{H_B}-u_{H_{\BB}}|=|u_{H_{\BB}}|\le \tg(x)
\le \frac{C}{\delta_S}(d(x,S)+d(x,y)+d(y,S))
(\tg(x)+\tg(y))
$$
proving \rf{difub}. In the same way we prove \rf{difub}
for the case $r_{\BB}>\delta_S,~r_B\le \delta_S$.
The remaining case $r_{\BB}>\delta_S,~r_B>\delta_S$
is trivial because here $u_{H_B}=u_{H_{\BB}}=0$.
\par We prove \rf{difuy} by a slight modification of the
proof given above. Using estimate \rf{fory} rather than
\rf{forall} we have
$$
|u_{H_B}-u(y)|\le C(d(x,S)+d(x,y))
(\tg(x)+\tg(y)).
$$
But $d(x,S)\le d(x,y)$, and \rf{difuy} follows.\bx
\begin{lemma}\label{fl}
Let $\BB\in W_S$ and let
$x\in \BB^*(:=\frac{9}{8}\BB)$. Then for every
$y\in X\setminus S$ we have
\bel{i1}
|\tu(x) - \tu(y)|\le
C\,\max_{B\in A}\min\{1, d(x,y)/r_{B}\}
|u_{H_B}-u_{H_{\BB}}|
\ee
where $A:=\{B\in W_S:~B^*\cap\{x,y\}\ne\emptyset\}$.
\par If $y\in S$, then
\bel{Iny}
|\tu(x)-u(y)|\le
C\max\{|u_{H_B}-u(y)|:B\in W_S, B^*\ni x\}.
\ee
\end{lemma}
{\bf  Proof.} By definition \rf{Extu} and properties
of the partition of unity we have
\be
I:=|\tu(x) - \tu(y)|&=&
\left|\sum_{B\in W_S} u_{H_B}\varphi_B(x)
-\sum_{B\in W_S} u_{H_B}\varphi_B(y)\right|\nn\\
&=&
\left|\sum_{B\in W_S} (u_{H_B}-u_{H_{\BB}})
(\varphi_B(x)-\varphi_B(y))\right|
\nn\\
&\le&
\sum_{B\in A}|u_{H_B}-u_{H_{\BB}}|
|\varphi_B(x)-\varphi_B(y)|
\nn
\ee
so that by property (3) of Lemma \ref{Wadd} for every
$y\in X\setminus S$
\bel{inq}
I\le
2N\max_{B\in A}|u_{H_B}-u_{H_{\BB}}|
|\varphi_B(x)-\varphi_B(y)|.
\ee
Since $0\le\varphi_B\le 1$, this implies
$$
I\le C\max\{|u_{H_B}-u_{H_{\BB}}|:B\in A\}.
$$
On the other hand, by property (d) of partition of unity we
have
$$
I\le
C_1d(x,y)\max_{B\in A}r_B^{-1}
|u_{H_B}-u_{H_{\BB}}|.
$$
Clearly, these inequalities imply \rf{i1}.
Similarly to \rf{inq}, for $y\in S$ we have
$$
|\tu(x)-u(y)|\le N\max\{|u_{H_B}-u(y)|:B\in W_S, B^*\ni x\}
$$
proving \rf{Iny}.\bx
\par We are in a position to prove $C\tg$ for some $C$ is a
generalized gradient of $\tu$.
\begin{lemma}\label{UXY}
The inequality
$$
|\tu(x) - \tu(y)|\le Cd(x, y)(\tg(x)+\tg(y))
$$
holds $\mu$-a.e. on $X$.
\end{lemma}
{\bf  Proof.} We will suppose that $x,y\in S\setminus E$
where $E$ is a subset of $S$ from inequality \rf{lipug}
(recall that $\mu(E)=0$). Clearly, for $x,y\in S$
the result follows from \rf{lipug}
so we may assume that $x\in X\setminus S$.
We let $\BB\in W_S$ denote
a Whitney ball such that $\BB\ni x$.
\par Denote $I:=|\tu(x) - \tu(y)|$ and consider two cases.
\par {\it The first case: $y\in \BB^*$.}
Since $x\in \BB$, we have $d(x,y)\le 2r_{\BB^*}\le
3r_{\BB}.$ Moreover, by \rf{drx} $r_{\BB}\approx
d(x,S)\approx d(y,S)$ and by inequality \rf{i1}
$$
I\le Cd(x,y)
\max\{r_B^{-1}|u_{H_B}-u_{H_{\BB}}|:~B\in W_S,
B^*\cap\{x,y\}\ne\emptyset\}.
$$
Since $x,y\in \BB^*$, for every ball $B\in W_S$ such that
$B^*\cap\{x,y\}\ne\emptyset$ we have
$B^*\cap \BB^*\ne\emptyset$. Therefore by
\rf{eqvball} $r_{\BB}\approx r_B.$
In addition, by Lemma \ref{uBB}
$$
|u_{H_B}-u_{H_{\BB}}|\le C(d(x,S)+d(x,y)+d(y,S))
(\tg(x)+\tg(y))
$$
so that
$$
|u_{H_B}-u_{H_{\BB}}|\le C r_{\BB}(\tg(x)+\tg(y)).
$$
Hence
\be
I&\le&
Cd(x,y)r_{\BB}^{-1}
\max\{|u_{H_B}-u_{H_{\BB}}|:~B\in W_S,
B^*\cap\{x,y\}\ne\emptyset\}\nn\\
&\le& C d(x,y)(\tg(x)+\tg(y)).\nn
\ee
\par {\it The second case: $y\notin \BB^*$.} Since
$x\in \BB$, this implies
$d(x,y)\ge \frac{1}{8}r_{\BB}$. Recall that
$r_{\BB}\approx d(x,S)$ so that
$d(x,S)\le Cd(x,y)$.
Since the distance function $d(\cdot,S)$
satisfies the Lipschitz condition, we have
$$
d(y,S)\le d(x,S)+d(x,y)\le Cd(x,y).
$$
\par Let $y\notin S$. Then by \rf{i1}
$$
I\le
C\max\{|u_{H_B}-u_{H_{\BB}}|:~B\in W_S,
B^*\cap\{x,y\}\ne\emptyset\}
$$
so that by \rf{difub}
$$
I\le C(d(x,S)+d(x,y)+d(y,S))
(\tg(x)+\tg(y))\le
Cd(x,y)(\tg(x)+\tg(y)).
$$
\par In the remaining case, i.e., for  $y\in S$,
the lemma follows from estimates \rf{difuy} and
\rf{Iny}.\bx
\par Let $f\in L^p(S), 1\le p\le\infty$.
We define an extension $F$ of $f$ by letting
$F(x):=f(x),~x\in S,$ and
\bel{sum}
F(x):=\sum_{B\in W_S}\limits |f_{H_B}|\chi_{B^*},
~~~~x\in X\setminus S.
\ee
\begin{lemma}\label{LPH}
$\|F\|_{L^p(X)}\le C\|f\|_{L^p(S)}.$
\end{lemma}
{\bf  Proof.} We will prove the lemma for the case $1\le
p<\infty$; corresponding changes for $p=\infty$ are
obvious. By property (3) of Lemma \ref{Wadd} for every
$x\in X\setminus S$ at most $N=N(C_d)$ terms of the sum in
\rf{sum} are not equal zero. Therefore
$$
|F(x)|^p\le C\sum_{B\in W_S}\limits |f_{H_B}|^p
\chi_{B^*}(x),
~~~x\in X\setminus S.
$$
This inequality and the doubling condition imply
$$
\int_{X\setminus S}|F|^p d\mu
\le C\sum_{B\in W_S}\limits
|f_{H_B}|^p\mu(B^*)
\le
C\sum_{B\in W_S}\limits |f_{H_B}|^p\mu(B).
$$
Recall that $\mu(H_B)\approx \mu(B)$ whenever
$r_B\le\delta$, see (i), (ii), Theorem \ref{HB},
so that
$$
|f_{H_B}|^p=\left|\frac{1}{\mu(H_B)}
\int_{H_B}fd\mu\right|^p
\le
\frac{1}{\mu(H_B)}\int_{H_B}|f|^pd\mu
\le
C\frac{1}{\mu(B)}\int_{H_B}|f|^pd\mu.
$$
Recall also that $H_B=\emptyset$ if $r_B>\delta$. Hence
$$
\int_{X\setminus S}|F|^p d\mu
\le C\sum_{B\in W_S}\limits
\int_{H_B}|f|^p d\mu
=C\int_{S}|f|^p
\left(\sum_{B\in W_S}\limits\chi_{H_B}\right)d\mu
$$
so that by property (iii) of Theorem \ref{HB}
$$
\int_{X\setminus S}|F|^p d\mu\le
C\int_{S}|f|^pd\mu.
$$
It remains to note that $F|_S=f$ and the lemma follows.\bx
\par Let us prove that
\bel{LPN}
\|\tu\|_{L^p(X)}\le C\|u\|_{L^p(S)}.
\ee
Since $0\le\varphi_B\le 1$ for every $B\in W_S,$
and $\supp \varphi_B\subset B^*,$ by \rf{Extu}
for every $x\in X\setminus S$ we have
$$
|\tu(x)|=
|\sum_{B\in W_S}\limits u_{H_B}\varphi_{B}(x)|
\le \sum_{B\in W_S}\limits|u_{H_B}|\varphi_{B}(x)
\le \sum_{B\in W_S}\limits|u_{H_B}|\chi_{B^*}(x).
$$
Hence $|\tu|\le |F|$ where $F(x):=u(x)$ for $x\in S$ and
$$
F(x):=\sum_{B\in W_S}\limits|u_{H_B}|
\chi_{B^*}(x),~~~x\in X\setminus S.
$$
Thus $\|\tu\|_{L^p(X)}\le\|F\|_{L^p(X)}$. But by Lemma
\ref{LPH} $\|F\|_{L^p(X)}\le C\|u\|_{L^p(S)}$, and \rf{LPN}
follows.
\par It remains to estimate $L^p$-norm of $\tg$.
To this end we define a function $G$ by letting
$G(x):=g(x),~x\in S$ and
$$
G(x):=\sum_{B\in W_S}\limits|g_{H_B}|\chi_{B^*}(x),
~~~x\in  X\setminus S.
$$
Then by \rf{defG} $|\tg|\le |G|+|F|$. By Lemma \ref{LPH}
$\|G\|_{L^p(X)}\le C\|g\|_{L^p(S)}$ so that
$$
\|\tg\|_{L^p(X)}\le \|G\|_{L^p(X)}+\|F\|_{L^p(X)}
\le C(\|g\|_{L^p(S)}+\|u\|_{L^p(S)}).
$$
\par Theorem \ref{EXT2} is completely proved.
\SECT{4. The sharp maximal function: proof of Theorems
\ref{EXT1} and \ref{EXT3}.}{4}
\indent
\par Let us fix a ball $K=B(z,r)$ such that
$K\cap S\ne\emptyset$. We denote
two families of balls associated to $K$
by letting
$
\BK:=\{B\in W_S:~B^*\cap K\ne\emptyset\}
$
and
$$
\BT:=\{B\in W_S:~B^*\cap K\ne\emptyset,
~r_B\le \delta_S\}.
$$
\begin{lemma}\label{e0} (i). For every $c\in\R$
$$
\int_{K\setminus S}|\tu-c|\,d\mu \le C\sum_{B\in
\BK}\mu(B)|u_{H_B}-c|.
$$
\par (ii). For every ball $B\in \BK$ we have
$r_B\le \eta_1 r$.
\par (iii). For every $c\in\R$
$$
\sum_{B\in \BT}\mu(B)|u_{H_B}-c|\le C\int_{(\eta_2K)\cap
S}|u-c|\,d\mu.
$$
\par Here $\eta_1,\eta_2$ are constants depending
only on the doubling constant $C_d$.
\end{lemma}
{\bf  Proof.} Let us prove property (i).
Recall that
$\sum\{\varphi_B(x):~B\in W_S\}=1$ for every
$x\in X\setminus S$. Then by definition \rf{Extu}
\be I:=\int_{K\setminus S}|\tu-c|\,d\mu &=&
\int_{K\setminus S}
|\sum_{B\in W_S}u_{H_B}\varphi_B-c|\,d\mu\nn\\
&\le&
\int_{K\setminus S}\sum_{B\in W_S}
|u_{H_B}-c|\varphi_B d\mu
=
\sum_{B\in W_S}\int_{K\setminus S}
|u_{H_B}-c|\varphi_B d\mu.\nn
\ee
Hence, by properties (a),(b) of the partition of unity and
by the doubling condition
$$
I\le
\sum_{B\in \Bc_K}\int_{B^*}
|u_{H_B}-c|\varphi_B d\mu
\le
\sum_{B\in \Bc_K}\mu(B^*)|u_{H_B}-c|
\le
C\sum_{B\in \Bc_K}\mu(B)|u_{H_B}-c|.
$$
\par Prove (ii). Let $B\in \Bc_K$ and let
$y\in B^*\cap K$. Then by \rf{drx}
$B^*\subset X\setminus S$ so that $y\notin S$.
Therefore there is a ball $B'\in W_S$
which contains $y$. Since $K\cap S\ne\emptyset$ and
$B'\cap K\ne \emptyset$, we have $\dist(B',S)\le 2r$.
But by Theorem \ref{Wcov}
$r_{B'}\le \dist(B',S)$ so that $r_{B'}\le 2r$.
In addition, $(B')^*\cap B^*\ne\emptyset$
so that by \rf{eqvball} $r_{B'}\approx r_B$.
This implies the required inequality $r_B\le \eta_1 r$
with some constant $\eta_1=\eta_1(C_d)$.
\par Prove (iii). We denote $A:=\cup\{H_B:~B\in\BT\}$ and
$$
m_K(x):=\sum\{\chi_{H_B}(x):~B\in \BT\}.
$$
Since $|u_{H_B}-c|\le |u-c|_{H_B}$ and
$\mu(H_B)\approx \mu(B)$, see (ii), Theorem \ref{HB},
\be \sum_{B\in \BT}\mu(B)|u_{H_B}-c| &\le& \sum_{B\in
\BT}\frac{\mu(B)}{\mu(H_B)} \int_{H_B}|u-c|\,d\mu \le
\gamma_2\sum_{B\in \BT}~~
\int_{H_B}|u-c|\,d\mu\nn\\
&=&
\gamma_2\int_{A}|u-c|m_Kd\mu.\nn
\ee
By property (i) of Theorem \ref{HB} for every $B\in\BT$
we have $H_B\subset(\gamma_1B)\cap S$.
Since $B^*\cap K\ne\emptyset$ and $r_B\le \eta_1 r$,
we obtain
$$
(\gamma_1B)\subset (1+(\gamma_1+9/8)\eta_1)K=\eta_2 K
$$
so that $H_B\subset (\eta_2 K)\cap S$. Thus
$A\subset (\eta_2 K)\cap S$.
\par It remains to note that by property (iii) of
Theorem \ref{HB} $m_K\le\gamma_3$ and the required
property (iii) follows.\bx
\begin{lemma}\label{e1} For every ball $K=B(z,r)$ such
that $z\in S$ and $r\le \delta_S/\eta_1$ we have
$$
\frac{r^{-\alpha}}{\mu(K)} \int_K |\tu-\tu_{K}|\,d\mu \le
Cu^\sharp_{\alpha,S}(z).
$$
\end{lemma}
{\bf  Proof.} We denote $D:=(\eta_2K)\cap S$ where $\eta_2$
is the constant from inequality (iii) of Lemma \ref{e0}.
Let us prove that
\bel{ei}
\int_K |\tu-\tu_K|\,d\mu \le C \int_D|u-u_D|\,d\mu
\ee
\par Since $r\le \delta_S/\eta_1$, by
(ii) of Lemma \ref{e0} we have
$r_B\le \delta_S$ for every ball $B\in \Bc_K$.
Thus $\BK=\BT$ so that $\{H_B:~B\in\BK\}$ is a subfamily of
the family $\Hc_S$ satisfying properties
(i)-(iii) of Theorem \ref{HB}.
\par Applying property (i) of Lemma \ref{e0}
with $c:=u_D$ we obtain
$$
\int_{K\setminus S}|\tu-c|\,d\mu \le C_d\sum_{B\in
\BT}\mu(B)|u_{H_B}-c|
$$
so that by (iii) of Lemma \ref{e0}
$$
\int_{K\setminus S}|\tu-c|\,d\mu\le C\int_D|u-c|\,d\mu.
$$
This implies
$$
\int_K |\tu-c|\,d\mu=\int_{K\cap S} |u-c|\,d\mu
+\int_{K\setminus S} |\tu-c|\,d\mu \le C\int_D|u-c|\,d\mu
$$
so that
$$
\int_K|\tu-\tu_{K}| \le 2\int_{K}|\tu-c|\,d\mu \le
C\int_D|u-c|\,d\mu
$$
proving \rf{ei}. Since
$\mu(K)\approx \mu(\eta_2K)$,
we finally obtain
$$
\frac{r^{-\alpha}}{\mu(K)} \int_K |\tu-\tu_{K}|\,d\mu \le
C(\eta_2r)^{-\alpha} \left(\frac{1}{\mu(\eta_2K)}
\int_D|u-u_D|\,d\mu\right) \le C u^\sharp_{\alpha,S}(z).\BX
$$
\par Recall that given a function $u$ defined on $S$ we let
$\br{u}$ denote its extension by $0$ to all of $X$.
As usual given $f\in L^{1,loc}(X)$ we let $\Mc f$ denote
the Hardy-Littlewood maximal operator:
$$
\Mc f(x):=\sup_{r>0}\frac{1}{\mu(B(x,r))}
\int_{B(x,r)}|f|\,d\mu.
$$
\begin{lemma}\label{e2}
Let $K=B(z,r)$ be a ball such that $z\in S$ and
$r>\delta_S/\eta_1$. Then
$$
\frac{r^{-\alpha}}{\mu(K)} \int_K |\tu-\tu_{K}|\,d\mu \le C
\Mc\br{u}(z).
$$
\end{lemma}
{\bf  Proof.} Applying property (i) of Lemma \ref{e0}
with $c:=0$ we obtain
$$
\int_{K\setminus S}|\tu|\,d\mu \le C\sum_{B\in
\BK}\mu(B)|u_{H_B}|.
$$
Since $u_{H_B}:=0$ whenever $r_B>\delta_S$, we have
$$
\int_{K\setminus S}|\tu|\,d\mu \le C\sum_{B\in
\BT}\mu(B)|u_{H_B}|.
$$
Applying (iii) of Lemma \ref{e0} with $c:=0$ we obtain
$$
\int_{K\setminus S}|\tu|\,d\mu \le C\int_{(\eta_2K)\cap
S}|u|\,d\mu.
$$
Since $r> \delta_S/\eta_1$, this implies
\be I:=\frac{r^{-\alpha}}{\mu(K)} \int_K
|\tu-\tu_{K}|\,d\mu &\le& 2
\frac{r^{-\alpha}}{\mu(K)}\int_K|\tu|\,d\mu \le
\frac{2\eta_1^\alpha}{\delta_S^\alpha}
\frac{1}{\mu(K)}\int_K|\tu|\,d\mu\nn\\
&\le& \frac{2\eta_1^\alpha}{\delta_S^\alpha}
\frac{1}{\mu(K)} \left(\int_{K\cap S}|u|\,d\mu+
C\int_{(\eta_2K)\cap S}|u|\,d\mu\right)\nn \ee
so that
$$
I\le \frac{4\eta_1^\alpha C}{\delta_S^\alpha}
\frac{1}{\mu(K)} \int_{(\eta_2K)\cap S}|u|\,d\mu.
$$
Since $\mu(K)\approx\mu(\eta_2K)$, we have
$$
I\le \frac{C}{\mu(\eta_2 K)} \int_{(\eta_2 K)\cap
S}|u|\,d\mu = \frac{C}{\mu(\eta_2 K)} \int_{\eta_2
K}|\br{u}|\,d\mu \le C \Mc\br{u}(z).\BX
$$
\par Lemma \ref{e1} and Lemma \ref{e2} imply
the following
\begin{proposition}\label{pr1}
For every $z\in S$
$$
(\tu)^\sharp_\alpha(z)\le
C(u^\sharp_{\alpha,S}(z)+\Mc\br{u}(z)).
$$
\end{proposition}
\par Let us estimate the value of
$(\tu)^\sharp_\alpha(z)$ for $z\in X\setminus S$. We will
put $\inf\limits_{H_Q}u^\sharp_{\alpha,S}:=0$ whenever
$H_Q=\emptyset$ (recall that $H_Q=\emptyset$ iff
$r_Q>\delta_S$).
\begin{lemma}\label{est1}
Let $Q=B(x_Q,r_Q)\in W_S$ and let $z\in Q$. Then
for every ball $K:=B(z,r)$ with $r\le\frac{1}{8}r_Q$
we have
\bel{fin}
\frac{r^{-\alpha}}{\mu(K)} \int_K
|\tu-\tu_{K}|\,d\mu \le
C\,\left(\inf_{H_Q}u^\sharp_{\alpha,S}
+\Mc\br{u}(z)\right).
\ee
\end{lemma}
{\bf  Proof.} We have to prove that for
arbitrary $s\in H_Q$
\bel{sum1}
I:=\frac{r^{-\alpha}}{\mu(K)} \int_K
|\tu-\tu_{K}|\,d\mu \le C\,(u^\sharp_{\alpha,S}(s)
+\Mc\br{u}(z)).
\ee
Since $r\le\frac{1}{8}r_Q$, the ball $K=B(z,r)\subset
\frac{9}{8}Q=:Q^*$. By Lemma \ref{fl} for every $x,y\in
K(\subset Q^*)$
$$
|\tu(x)-\tu(y)|\le
Cd(x,y)\max\{r_B^{-1}|u_{H_B}-u_{H_Q}|:~B\in W_S,
B^*\cap \{x,y\}\ne\emptyset\}.
$$
Since $x,y\in Q^*$, for every ball $B\in W_S$ such that
$B^*\cap \{x,y\}\ne\emptyset$ we have
$B^*\cap Q^*\ne\emptyset$. Therefore by \rf{eqvball}
\bel{BapQ}
\frac{1}{C_1}r_Q\le r_B\le C_1r_Q.
\ee
We denote $A:=\{B\in W_S:~B^*\cap Q^*\ne\emptyset\}$. Then
$$
|\tu(x)-\tu(y)|\le
C\frac{d(x,y)}{r_Q}
\max_{B\in A}|u_{H_B}-u_{H_Q}|.
$$
Hence
\be \frac{1}{\mu(K)}\int_K|\tu-\tu_K|\,d\mu &\le&
\frac{1}{\mu(K)^2}\int_K\int_K
|\tu(x)-\tu(y)|\,d\mu(x)d\mu(y)\nn\\
&\le&
C\frac{d(x,y)}{r_Q}
\max_{B\in A}|u_{H_B}-u_{H_Q}|.\nn
\ee
Since $d(x,y)\le \diam K\le 2r$ and
$r\le\frac{1}{8}r_Q$, this implies
\bel{uKB}
I\le Cr_Q^{-\alpha}\max_{B\in A}|u_{H_B}-u_{H_Q}|.
\ee
\par Let us consider two cases.
\par {\it The first case:} $r_Q\le\delta_S/C_1$ where $C_1$
is the constant from inequality \rf{BapQ}. Then for each
$B\in A$ we have $r_B\le\delta_S$ so that $H_B,H_Q$ satisfy
properties (i),(ii) of Theorem \ref{HB}. Thus
$H_B\subset(\gamma_1B)\cap S$,~ $H_Q\subset(\gamma_1Q)\cap
S$,~ $\mu(H_B)\approx \mu(B)$, and $\mu(H_Q)\approx
\mu(Q)$.
\par Since $B^*\cap Q^*\ne\emptyset$ and $r_B\approx r_Q$,
for some positive $C_2=C_2(\gamma_1)$ we have
$$
B\cup Q\cup H_B\cup H_Q\subset D:=B(s,C_2r_Q).
$$
(Recall that $s$ is an arbitrary point of $H_Q$.)
These inequalities and the doubling condition imply
$\mu(H_B)\approx \mu(H_Q)\approx \mu(D)$.
Hence
$$
|u_{H_B}-u_{D\cap S}| \le \frac{1}{\mu(H_B)}\int_{H_B}
|u-u_{D\cap S}|\,d\mu\le C\frac{1}{\mu(D)}\int_{D\cap S}
|u-u_{D\cap S}|\,d\mu.
$$
A similar estimate is true for $H_Q$ so that
$$
|u_{H_B}-u_{H_Q}|\le |u_{H_B}-u_{D\cap
S}|+|u_{H_Q}-u_{D\cap S}| \le C\frac{1}{\mu(D)}\int_{D\cap
S} |u-u_{D\cap S}|\,d\mu.
$$
Applying this inequality to \rf{uKB} we obtain
$$
I\le C\frac{r_Q^{-\alpha}}{\mu(D)}\int_{D\cap S}
|u-u_{D\cap S}|\,d\mu \le
C\frac{r_D^{-\alpha}}{\mu(D)}\int_{D\cap S} |u-u_{D\cap
S}|\,d\mu
$$
where $r_D:=C_2r_Q$ is the radius of the ball
$D:=B(s,C_2r_Q)$. Hence by definition \rf{smfS} we have
$I\le C u^\sharp_{\alpha,S}(s)$ proving \rf{sum1}.
\par {\it The second case:} $r_Q>\delta_S/C_1$.
By \rf{uKB} $I\le C\max\{|u_{H_B}|:~B\in A\}$.
Recall that $u_{H_B}:=0$ if $r_B>\delta_S$ so that
$$
I\le C\max\{|u_{H_B}|:~B\in A, r_B\le\delta_S\}.
$$
\par By Theorem \ref{HB} for every $B\in A$
such that $r_B\le\delta_S$ we have $H_B\subset
(\gamma_1B)\cap S$, ~$\mu(H_B)\approx \mu(B)$. Since
$r_B\approx r_Q$ and $z\in Q$, for some positive
$C_3=C_3(\gamma_1)$ we have $H_B\subset B(z,C_3r_Q)$. Put
$\widetilde{D}:=B(z,C_3r_Q)$. Since
$\mu(\widetilde{D})\approx\mu(Q)$ and
$\mu(B)\approx\mu(Q)$, we have
$\mu(H_B)\approx\mu(\widetilde{D})$. Hence
$$
|u_{H_B}|\le\frac{1}{\mu(H_B)} \int_{H_B}|u|\,d\mu\le
C\frac{1}{\mu(\widetilde{D})} \int_{\widetilde{D}\cap
S}|u|\,d\mu \le C\Mc\br{u}(z)
$$
proving that $I\le C\Mc\br{u}(z)$.\bx
\begin{lemma}\label{est2} Inequality \rf{fin} is
true for every $r>\frac{1}{8}r_Q$.
\end{lemma}
{\bf  Proof.} We denote $\eta_3:=8(\gamma_1+10)$,
$\tr:=\eta_3r$ and $\tK:=\eta_3K=B(z,\tr)$.
Recall that $\gamma_1$ is the constant
from Theorem \ref{HB}. Prove that
$\tK\cap S\ne\emptyset$.
In fact, let $a_Q\in Q$ and $b_Q\in S$
be points satisfying the inequality
$d(a_Q,b_Q)\le 2d(Q,S)$.
Then by (ii), Theorem \ref{Wcov},
$d(a_Q,b_Q)\le 2d(Q,S)\le 8 r_Q$. But $z\in Q$ so that
$$
d(z,b_Q)\le d(z,a_Q)+d(z,b_Q)\le 2r_Q+8 r_Q=10r_Q
\le 80r\le\eta_3 r=\tr.
$$
Thus $b_Q\in\tK\cap S$ proving that $\tK\cap
S\ne\emptyset$.
\par Let us consider two cases.
\par {\it The first case:} $\tr:=\eta_3r>\delta_S/\eta_1.$
Since $\tr$ is the radius of the ball $\tK=\eta_3K$,
$\tr>\delta_S/\eta_1$ and $\tK\cap S\ne\emptyset$,
by Lemma \ref{e2}
$$
\frac{\tr^{-\alpha}}{\mu(\tK)} \int_{\tK}
|\tu-\tu_{\tK}|\,d\mu \le C \Mc\br{u}(z).
$$
By the doubling condition $\mu(K)\approx \mu(\tK)$ so that
\be I:=\frac{r^{-\alpha}}{\mu(K)} \int_{K}
|\tu-\tu_{K}|\,d\mu &\le& \frac{2r^{-\alpha}}{\mu(K)}
\int_{K}
|\tu-\tu_{\tK}|\,d\mu\nn\\
&\le& C\frac{\tr^{-\alpha}}{\mu(\tK)} \int_{\tK}
|\tu-\tu_{\tK}|\,d\mu \le C\Mc\br{u}(z)\nn \ee
proving \rf{fin}.
\par {\it The second case.} $\tr:=\eta_3r\le\delta_S/\eta_1.$
Since $8r>r_Q$, we have
$$
r_Q<8\delta_S/(\eta_1\eta_3)<\delta_S.
$$
Therefore by Theorem \ref{HB} $H_Q\ne\emptyset$,
$\mu(H_Q)\approx \mu(Q)$ and
$H_Q\subset(\gamma_1Q)\cap S$.
\par Take $s\in H_Q$ and put $V:=B(s,\eta_3r)$.
Since $H_Q\subset\gamma_1Q$,
$d(s,x_Q)\le\gamma_1r_Q$ so that for every
$a\in K=B(z,r)$
\be
d(s,a)&\le&
d(s,x_Q)+d(x_Q,z)+d(z,a)\nn\\
&\le& \gamma_1r_Q+r_Q+r\le 8\gamma_1r+8r+r=
(8\gamma_1+9)r\le \eta_3r\nn
\ee
proving that $K\subset V$. On the other hand,
$V\subset 2\eta_3K$ so that by the doubling condition
$\mu(V)\approx \mu(K)$. Hence
$$
I\le 2r_K^{-\alpha}
\frac{1}{\mu(K)}\int_K|\tu-\tu_{V}|\,d\mu \le
Cr_V^{-\alpha} \frac{1}{\mu(V)}\int_V|\tu-\tu_{V}|\,d\mu.
$$
But $r_V:=\eta_3r\le \delta_S/\eta_1$ so that
by Lemma \ref{e1}  $I\le Cu^\sharp_{\alpha,S}(s)$.
This finishes the proof of \rf{fin} and the lemma.\bx
\begin{theorem}\label{pr}
For every $z\in X$
$$
(\tu)^\sharp_\alpha(z)\le
C(\Mc(u^\sharp_{\alpha,S})^\cw(z)+
\Mc u^\cw(z)).
$$
\end{theorem}
{\bf  Proof.} For $z\in S$ this follows
from Proposition \ref{pr1}.
\par Let $Q\in W_S$ and let $z\in Q$.
Then by Lemmas \ref{est1} and \ref{est2}
\bel{A1}
(\tu)^\sharp_\alpha(z)\le
C\,\left(\inf_{H_Q}u^\sharp_{\alpha,S}
+\Mc\br{u}(z)\right).
\ee
Recall that in this formula we put the infimum
to be equal $0$ whenever $H_Q=\emptyset$,
i.e., $r_Q>\delta_S$. Therefore in the
remaining part of the proof we may assume that
$r_Q\le \delta_S$. Then by Theorem \ref{HB}
$\mu(H_Q)\approx \mu(Q)$ and
$H_Q\subset(\gamma_1 Q)\cap S$.
\par Let us denote $B:=B(z,(\gamma_1+1)r_Q)$
and $h:=(u^\sharp_{\alpha,S})^\cw$. Since $z\in Q$,
we have $H_Q\subset \gamma_1 Q\subset B$. In addition,
by the doubling condition
$\mu(H_Q)\approx \mu(B).$ Hence
$$
\inf_{H_Q}u^\sharp_{\alpha,S}=
\inf_{H_Q}h
\le
\frac{1}{\mu(H_Q)}\int_{H_Q}hd\mu
\le\frac{1}{\mu(H_Q)}\int_{B}hd\mu\nn\\
\le \frac{C}{\mu(B)}\int_{B}hd\mu\le C\Mc h(z).
$$
This inequality and \rf{A1} imply the proposition.\bx
\begin{remark}
{\em Similar estimates and definition of $\tu$, see
\rf{Extu}, easily imply that $|\tu(x)|\le C\Mc u^\cw(x)$
for every $x\in X$.}
\end{remark}
\par {\bf Proof of Theorem \ref{EXT3}.} It can be easily
shown that for any extension $U$ of a function
$u\in L^p(S)$ to all of $X$ we have
$u_{\alpha,S}^{\sharp}(x)\le
2U_{\alpha}^{\sharp}(x),~ x\in S$. This  immediately
implies the inequality
$$
\|u\|_{L^p(S)}+\|u_{\alpha,S}^{\sharp}\|_{L^p(S)}
\le 2\|u\|_{\Cc^\alpha_p(X,d,\mu)|_S}.
$$
\par Now let $u,u^\sharp_{\alpha,S}\in L^p(S),
1<p\le\infty$. Prove that
$\tu=\Ext u\in\Cc^\alpha_p(X,d,\mu)$.
By Theorem \ref{pr}
$$
\|(\tu)^\sharp_\alpha\|_{L^p(X)}
\le
C(\|\Mc(u^\sharp_{\alpha,S})^\cw\|_{L^p(X)}+
\|\Mc u^\cw\|_{L^p(X)}).
$$
Recall that the operator $\Mc$ is bounded
in $L^p(X)$ whenever $1<p\le \infty$ and $(X,d,\mu)$ is
a metric space of a homogeneous type, see,
e.g. \cite{He}, p. 10. Hence
$$
\|(\tu)^\sharp_\alpha\|_{L^p(X)}\le
C(\|(u^\sharp_{\alpha,S})^\cw\|_{L^p(X)}+
\|u^\cw\|_{L^p(X)})
=C(\|u^\sharp_{\alpha,S}\|_{L^p(S)}+\|u\|_{L^p(S)}).
$$
Since $\|\tu\|_{L^p(X)}\le C\|u\|_{L^p(S)}$, see \rf{LPN},
we finally obtain
$$
\|\tu\|_{\Cc^\alpha_p(X,d,\mu)}:=
\|\tu\|_{L^p(X)}+\|(\tu)^\sharp_\alpha\|_{L^p(X)} \le
C(\|u\|_{L^p(S)}+\|u^\sharp_{\alpha,S}\|_{L^p(S)})
$$
proving that $\tu\in \Cc^\alpha_p(X,d,\mu)$ and equivalence
\rf{Eq} holds.
\par The proof of Theorem \ref{EXT3} is complete.\bx

\end{document}